\documentclass[12pt]{article}

\usepackage{amssymb}
\usepackage{amscd}
\usepackage{epsfig}
\usepackage{graphicx}
\usepackage{theorem}
\usepackage{amsmath}

\input xy
\xyoption{all}

\newtheorem{thm}{Theorem}[section]

\newtheorem{prop}[thm]{Proposition}

\newtheorem{lemm}[thm]{Lemma}
\newtheorem{conj}[thm]{Conjecture}

{\theorembodyfont{\rm}

\newtheorem{defn}[thm]{Definition}
\newtheorem{exam}[thm]{Example}
\newtheorem{rem}[thm]{Remark}

}

\newenvironment{proof}{\noindent {\em Proof.}}{\hfill 
$\square$\newline}

\newcommand{\no}{\noindent}

\newcommand{\La}{\Lambda}

\newcommand{\al}{\alpha}

\newcommand{\Q}{\mathbb Q}
\newcommand{\Z}{\mathbb Z}
\newcommand{\R}{\mathbb R}
\newcommand{\fO}{\mathfrak o}

\def\cX{{\mathcal X}}

\def\cC{{\mathcal C}}
\def\cD{{\mathcal D}}

\def\cH{{\mathcal H}}

\def\cO{{\mathcal O}}

\def\cL{{\mathcal L}}
\def\cM{{\mathcal M}}
\def\cN{{\mathcal N}}

\def\cU{{\mathcal U}}

\def\oM{{\overline{\mathcal M}}}

\newcommand{\PGL}{{\rm PGL}} 
\newcommand{\SL}{{\rm SL}}

\def\Sym{{\rm Sym}}
\def\ovl{\overline}

\def\ra{\rightarrow}

\def\P{{\mathbb P}}
\def\Q{{\mathbb Q}}

\def\Z{{\mathbb Z}}

\def\N{{\mathbb N}}

\def\bP{{\mathbb P}}
\def\bQ{{\mathbb Q}}

\def\fS{{\mathfrak S}}

\def\Pic{{\rm Pic}}

\def\eps{{\epsilon}}

\makeatother
\makeatletter

\ifx\pdfoutput\undefined\else
\usepackage[pdftex,verbose,
           colorlinks=true,
             pdftitle={Integral points and effective cones of moduli spaces of stable maps},
             pdfauthor={Br. Hassett and Yu. Tschinkel}
            ]
              {hyperref}
 \fi

\author{
Brendan Hassett 
\thanks{Partially supported by NSF grant 0196187}
and 
Yuri Tschinkel 
\thanks{Partially supported by NSF grant 0100277 
and the Clay Foundation}}

\title{Integral points and effective cones of \\
moduli spaces of stable maps}

\begin{document} 
 
\date{May 8, 2002}



\maketitle

\begin{abstract}
Consider the Fulton-MacPherson configuration space
of $n$ points on $\P^1$, which is 
isomorphic to a certain moduli space of stable maps to $\P^1$.
We compute the cone of effective
${\mathfrak S}_n$-invariant divisors on this space.
This yields a geometric interpretation of known asymptotic 
formulas for the number of integral points of
bounded height on compactifications of $\SL_2$ 
in the space of binary forms of degree $n\ge 3$.  
\end{abstract}

\tableofcontents

\setcounter{section}{-1}

\section{Introduction}
\label{sect:introduction}

In this paper, we compute the ${\mathfrak S}_n$-invariant 
cone of effective divisors of the Fulton-MacPherson
configuration space of $n$ points on $\P^1$.  This space
is isomorphic to the moduli space $\oM_{0,n}(\P^1,1)$
of stable maps of degree one from genus zero curves with 
$n$ marked points to $\P^1$. 
We also compute the effective cone of the generic fiber of the
natural map 
$$
\oM_{0,n}(\P^1,1)/{\mathfrak S}_n \ra \oM_{0,n}/{\mathfrak S}_n.
$$
Our motivation is to provide a geometric explanation 
of a formula, obtained by Duke, Rudnick and Sarnak, 
giving the asymptotic behavior of the number of binary forms
of degree $n$ with fixed discriminant and bounded integral coefficients.
This fits into a larger program to predict and prove
asymptotic formulas for the number of rational and integral points
of bounded height on algebraic varieties.

\

We introduce a counting function for integral points on an algebraic variety
as follows: given a variety $U$ over a ring of integers 
$\fO$ and functions $g_1,...,g_n$, regular
on $U$, define 
$$
N(U,B):=\{ x\in U({\fO})\,|\, \max_j (\|g_j(x)\|)\le B\},
$$
where $\|\cdot \|$ is a valuation on $\fO$. 
This is finite only when the functions $g_j$ give an 
embedding of $U$.

It is most natural to interpret the functions $g_j$ 
as {\em sections} of a line bundle $L$ on a projective
compactification $X\supset U$  defined over the 
fraction field $F$ of $\fO$. 
The fact that the sections embed $U$ implies that 
$L$ is {\em big}, i.e., is contained in the interior of the 
effective cone $\La_{\rm eff}(X)$ of $X$. 
Therefore, in order to describe all natural counting functions
on open subsets of $X$ we need to compute its effective cone. 
Furthermore, in many cases it can be proved that 
the asymptotic properties of $N(U,B)$ 
are intimately related to the structure of this cone.

\

Let $P({\bf x})=P(x_0,...,x_r)$ be a homogeneous polynomial of degree $n$ in
$r+1$ variables. 
A standard heuristic in number theory predicts that 
the number 
$$
N_P(B):=\left\{ {\bf x}\,|\, \max (|x_j|)\le B, \,\, 
P({\bf x})=0,\,\, {\rm and}\,\, {\bf x}\in \Z^{r+1}\right\}
$$
of integral solutions  of the equation $P({\bf x})=0$ 
of ``height'' $\le B$ grows asymptotically like $B^{r+1-n}$
as $B\ra \infty$. When the number of variables is 
$\gg 2^n$, the affine variety $V_P$ defined by $P=0$
is smooth and there are no local obstructions,
an asymptotic formula can be established using
the classical circle method in analytic number theory 
(see \cite{Birch}, \cite{Schmidt} and the references therein).
Of course, there may be difficulties when the number of variables is small
or the variety $V_P$ is singular.

\

The following example appeared in the paper  
by Duke, Rudnick and Sarnak \cite{duke}.
Consider the vector space of binary forms of degree $n$
$$x_nz^n+x_{n-1}z^{n-1}w+\ldots + x_0 w^n.$$
The algebraic group $\SL_2$ acts on this space by coordinate
substitutions.
When $n=3$, the {\em discriminant form}
$$
{\rm disc}(x_0,...,x_3)
:=27x_0^2x_3^2-18x_0x_1x_2x_3+4x_0x_2^3+4x_1^3x_3-x_1^2x_2^2
$$
generates the ring of $\SL_2$-invariants.  
Then there exists a constant $c>0$ so that
$$
N_{{\rm disc}-1}(B)= c B^{2/3}(1+o(1))
$$
as $B\ra \infty$.  Note that the exponent $2/3$
is larger than what is predicted by the standard heuristic.  

More generally, one has the 
\begin{thm}\label{thm:example} \cite{duke}
Fix a generic binary form $f$ of degree $n\ge 3$ with integral coefficients.
Let $N(B)$ be the number of binary forms $\SL_2(\Z)$-equivalent to $f$ 
with coefficients bounded by $B$. 
Then there exists a $c>0$ such that
$$
N(B)=cB^{2/n}(1+o(1)),
$$
as $B\ra \infty$. 
\end{thm}

We give a geometric interpretation of the exponent $2/n$ in 
Theorem \ref{thm:example}.  To this end, we refine the
heuristics for counting integral points to take into account
singularities of the relevant varieties (see Conjecture \ref{conj:a}).
We verify that Conjecture \ref{conj:a}
is consistent with Theorem \ref{thm:example} in Theorem \ref{thm:main}.
Its proof
involves the computations of effective cones alluded to above.  

\no {\bf Acknowledgments:} We are grateful to the 
Alfr\'ed R\'enyi Institute of the Hungarian Academy
of Sciences for organizing the conference at which much
of this work was done.

\section{Generalities}

\subsection{Singularities of pairs and effective cones}
\label{sect:gen}

We work over a field of characteristic zero.  
Let $X$ be a normal projective
variety with canonical class $K_X$
and let $D$ be a reduced effective Weil divisor of $X$.

\begin{defn}
\label{defn:good}
A {\em good} pair $(X,D)$ consists of 
a smooth projective variety $X$
and a strict normal crossings divisor $D$ in $X$. 
This means that all irreducible components of $D$ are smooth and
intersect transversally. 
\end{defn}

\

Let $(X,D)$ be a good pair and let
$\La_{\rm eff}(X)$ denote the closed
cone of effective divisors classes of $X$;  a divisor is
big exactly when its class is in the interior of this
cone.  
Define 
$$
a(L,D):=\inf\{ a\in \R\,|\, aL+(K_X+D)\in \La_{\rm eff}(X)\},
$$
where we identify line bundles and their divisor classes.
Note that $a(L,D)$ is a positive real number whenever $-(K_X+D)$
is big.  The constant $-a(L,D)$ is called the {\em log-Kodaira energy} 
of $L$ (see \cite{Fu}).

\

If $(X,D)$ is not good then resolution of singularities
implies the existence of a {\em good resolution}
$\rho\,:\, (\tilde{X},\tilde{D})\ra (X,D)$.  Precisely,
$(\tilde{X},\tilde{D})$ is a good pair, $\rho$
a birational projective morphism, and $\tilde{D}$ is the union of
the exceptional divisors of $\rho$ and the proper transform of $D$.
Recall that $(X,D)$ 
is {\em log-canonical} if $K_X+D$ is $\bQ$-Cartier and
$$
K_{\tilde{X}} +\tilde{D} \equiv_{\bQ} \rho^*(K_X+D) + \sum d_jE_j,
$$
where the $E_j$ are the exceptional divisors of $\rho$ and  
$d_j\ge 0$ for all $j$.

\begin{exam}
\label{exam:lc}
When $X$ is a smooth surface, $(X,D)$ is log-canonical
only when the curve $D$ is smooth or nodal.  
If $X$ is smooth of arbitrary dimension, 
$D$ must have at worse nodes in codimension one.
\end{exam}

If $L$ is a line bundle on $X$
put 
$$
a(L,D):=a(\rho^*L,D^t),
$$
where $D^t\subset \tilde{X}$ is the total transform of $D$.  
Note that $a(L,D)$ is computed on $\tilde{X}$.

\begin{prop}
\label{prop:ind}
Let $(X,D)$ be a log-canonical pair
and assume that $X-D$ has
canonical singularities.  If $L$ is a big line bundle on $X$
then
$$
a(L,D)=\inf\{ a\in \R\,|\, aL+(K_X+D)\in \La_{\rm eff}(X)\}.
$$ 
In particular, $a(L,D)$ does not depend on the choice of  
a desingularization.
\end{prop}

\begin{proof}
Choose a good resolution $\rho\,:\, (\tilde{X},\tilde{D})\ra (X,D)$, so that
$$
K_{\tilde{X}} + \tilde{D} -\sum d_jE_j=\rho^*(K_X+D)
$$
where $d_j\ge 0$, and $d_j\ge 1$ if $\rho(E_j)\not \subset D$.  
In particular, each exceptional
divisor not contained in the total transform $D^t$ has
log discrepancy $\ge 1$.  Therefore, we have
$$
\al \rho^*(L) +K_{\tilde{X}} + D^t -\sum d'_jE_j=\rho^*(K_X+D+\al L),
$$
with each $d'_j\ge 0$.  
For any $\Q$-Cartier divisor $M$ on $X$ and effective divisor
$\sum d'_jE_j$ supported in the exceptional locus of $\rho$, 
$M$ is effective iff $\rho^*(M)+\sum d'_jE_j$ is effective. 
\end{proof}

\begin{prop}
\label{prop:eff}
Let $(X_1,D_1)$ and $(X_2,D_2)$ be log-canonical pairs,
so that $X_1-D_1$ and $X_2-D_2$ have canonical singularities.  
Assume that 
$\pi\,:\, X_1\ra X_2$ is a finite dominant morphism so that
$$
\pi^*(K_{X_2}+D_2)=K_{X_1}+D_1.
$$
Let $L$ be a big divisor on $X_2$. 
Then $a(L,D_2)=a(\pi^*(L),D_1)$.  
\end{prop}
In fact, it suffices to assume that either $(X_1,D_1)$
or $(X_2,D_2)$ satisfies the singularity condition 
\cite[Sect.~20.3]{flab}.

\begin{proof}
Given a finite dominant morphism
$\pi:X_1\ra X_2$ and a $\Q$-Cartier divisor $M$ on $X_2$,
$M$ is effective iff $\pi^*(M)$ is effective.  Indeed,
the divisor $\pi_*\pi^*M$ is defined and equal to $\deg(\pi)M$.  
Combining this with Proposition \ref{prop:ind} gives the result.  
\end{proof}

\begin{rem}
Let $(X,D)$ be a log terminal pair so that
$X-D$ has singularities which are {\em not} canonical.
Then our definition of the Kodaira energy differs
slightly from Fujita's \cite{Fu}.  In applications
to integral points, we are interested in invariants of
the open variety $X-D$.  In Fujita's definition, on
passing from $(X,D)$ to a good resolution, any
exceptional divisors over $X-D$ with negative discrepancy 
must be added to the boundary.  This changes
the open variety.  
\end{rem}

\subsection{Integral points}

Retain the notation from the previous section and 
assume that $X$ and $D$ are defined over a number field $F$. 
Let $\fO_S$ denote the ring of integers of $F$, 
where $S$ is a finite set of nonarchimedean places of $F$. 
Fix models $\cX$ and $\cD$ flat and proper over 
the ring of integers $\fO_S$.
A $(D,S)$-integral point is an $\fO_S$-point of $(\cX - \cD)$.
In particular, if $D=\emptyset$ an integral point is the same
as a rational point on $X$.

Let $\cL$ be a very ample metrized line bundle on $X$, 
$U\subset X$ a Zariski open subset and $\cU$ a model of $U$ over $\fO_S$.  
Let $S$ be a finite set of places in $F$, including the archimedean places. 
Let 
$$
N(\cU,\cL,B):=\# \{ x\in \cU(\fO_S)\,|\, H_{\cL}(x)\le B\}
$$
denote the number of $(D,S)$ integral points on 
$U$ of $\cL$-height bounded by $B$.
A natural extrapolation of Vojta's conjecture about 
integral and rational points on 
varieties of (log-)general type  \cite{vojta} and 
Batyrev-Manin conjectures about rational points of bounded height
on Fano varieties \cite{FMT,BM} would be:

\begin{conj}
\label{conj:a}
For any $\eps>0$, there exists
a dense Zariski open subset $U\subset X$ such that
$$
N(\cU,\cL,B)\ll B^{a(L,D)+\epsilon}
$$
as $B\ra \infty$. If $-(K_X+D)$ is big then 
$$
N(\cU,\cL,B)\gg B^{a(L,D)-\epsilon},
$$ 
as $B\ra \infty$,
at least after a suitable finite extension of $F$ and $S$. 
\end{conj}

The statement is independent of the choice of $S$ and 
the choice of a metrization on $L$.

\

Many precise results about asymptotics of rational and integral 
points are currently available
(see, for example, \cite{FMT,BT,CLT,peyre,peyre-2,duke,eskin-1,eskin-2}
and the references therein). 
As far as we know, Conjecture~\ref{conj:a} is compatible with all 
of them. However, to actually check this compatibility one has
to compute the geometric invariants of (some resolution of) 
the pair $(X,D)$.
In particular, one has to determine the effective cone. 
This can be a formidable task even for rational varieties, e.g., 
like the moduli space of pointed rational curves 
$\ovl{\cM}_{0,n}$ (see \cite{ht}). 

\subsection{Computing effective cones}
\label{sect:compeff}

Let $X$ be a nonsingular projective variety, perhaps with
an action by a finite group $G$.  
We review strategies for computing the $G$-invariant effective 
cone $\La_{\rm eff}(X)^G$ and thus the effective cone
of the quotient $X/G$ (cf. \cite{KeMc}).  

A curve class $[C] \in N_1(X)$ is said to be {\em nef} if $[C].D\geq 0$
for each $D\in \La_{\rm eff}(X)$.  
A family of curves
passing through the generic point of $X$ is automatically
nef.  Indeed, consider a family $\cC \ra B$ of integral projective curves
in $X$ and an irreducible codimension-one subvariety
$D\subset X$.  If, for generic $b\in B$, the fiber $C_b \not \subset D$,
we have $[C_b].D\geq 0$.  

Fix a collection of effective divisors
$$\Gamma=\{A_1,\ldots,A_m\}$$ 
which we expect to generate $\La_{\rm eff}(X)^G$.  
To prove that $\Gamma$
generates the ($G$-invariant)
effective cone, it suffices to find a collection
of nef ($G$-invariant) curve classes
$$\Xi=\{C_1,\ldots,C_{\ell}\}$$
so that the cone generated by $\Gamma$ contains the dual
to the cone generated by $\Xi$.  

\

In section \ref{sect:fm}, we shall
use a refinement of this method (see \cite{CR}, \cite{Ru}).
A divisor $D\in \La_{\rm eff}(X)$
is {\em moving relative to $\Gamma$} if some multiple
of $D$ contains no element of $\Gamma$ as a fixed component.
Every effective divisor is a sum 
$$M+\sum_{i=1}^m A_id_i, \quad d_i\ge 0$$
where $M$ is moving relative to $\Gamma$.  To prove that
$\Gamma$ generates the effective cone, it suffices to 
show that $M$ is an effective sum of the $A_i$.  

A curve class is {\em nef relative to $\Gamma$} if $[C].M\geq 0$
for each $M$ which is moving relative to $\Gamma$.  Any
family of curves passing through the generic point of 
{\em some} $A_i$ is nef relative to $\Gamma$.
Consequently, to show that $\Gamma$ generates the effective
cone, it suffices to find a collection $\Xi$
of curve classes, nef relative to $\Gamma$,
so that the cone generated by $\Gamma$ contains the dual
to the cone generated by $\Xi$.  

\section{Construction of resolutions}
\subsection{Binary forms and $\SL_2$-orbit closures}

Let $V$ be a two-dimensional vector space with coordinates $z$ and $w$,
equipped with the standard $\SL_2$-action.
Let $\Sym^nV^*$ be the space of binary forms of degree $n$ 
$$
f=x_0z^n+x_1 z^{n-1}w+\ldots + x_n w^n.
$$
It carries an induced action of $\SL_2$ by substitution.

Associating to each form $f\ne 0$ its roots $\alpha_1,\ldots,
\alpha_n$ yields a map
$$(\Sym^nV^* - 0) \ra \bP(V)^n/{\mathfrak S}_n$$
and an identification $\bP(\Sym^nV^*)\simeq \bP(V)^n/{\mathfrak S}_n$.
The {\em discriminant} of a polynomial $f$ is a homogeneous
form in its coefficients $x_0,\ldots,x_n$ and defines a
divisor $D\subset X=\bP(\Sym^nV^*)$.

Now we may state our main result:
\begin{thm}[Computation of Kodaira Energy] \label{thm:main}
Let $f$ be a generic bilinear form of degree $n$,
$X_f\subset \bP(\Sym^nV^*)$
the closure of the $\SL_2$-orbit
through $f$, $D_f$ the intersection of the discriminant
with $X_f$, and $L$ the restriction of the standard
polarization to $X_f$.  Then we have
$$a(L,D_f)=2/n.$$
In particular,
Conjecture \ref{conj:a} is consistent with 
Theorem \ref{thm:example}.
\end{thm}
To prove this, we require a resolution
(i.e., a partial desingularization)
of $(X_f,D_f)$ on which we may evaluate 
$a(L,D_f)$ using Proposition \ref{prop:ind}.  
This resolution will be induced by 
a natural resolution of $(X,D)$.  

\begin{rem}
Example \ref{exam:lc} shows that $(X,D)$
is far from being log-canonical. 
When $n=3$, the discriminant
has cusps in codimension one:
a transverse slice 
$$z^3+bzw^2+cw^3$$
intersects the discriminant in the cuspidal curve
$$4b^3+27c^2=0.$$
\end{rem}

Our resolution of $(X,D)$ will be a ${\mathfrak S}_n$-quotient
of a natural desingularization for $(\bP(V)^n,\Delta)$,
where $\Delta$ is the diagonal, i.e., the points lying over the
discriminant.  Both admit interpretations as moduli spaces
of stable maps.

\subsection{Moduli spaces}
\label{sect:mod}

Fix an integer $n\ge 3$.  Let $\ovl{\cM}_{0,n}$ denote
the Knudsen-Mumford moduli space of stable curves of genus zero with
$n$ marked points \cite{Kn}
$$(C,p_1,\ldots,p_n).$$
Let $\ovl{\cM}_{0,n}(\bP^1,1)$ denote the 
Kontsevich moduli space
of stable maps of degree one from
genus-zero curves with $n$ marked points to $\bP^1$ \cite{K,KM,fp}
$$(C,p_1,\ldots,p_n,\mu:C \ra \bP^1).$$
This is naturally isomorphic to 
the Fulton-MacPherson \cite{fm}
configuration space $\bP^1[n]$ for $n$ points in $\bP^1$
(see \cite{fp} \S 0).  
However, for our purposes it is convenient to use the 
moduli space notation.

We have the following natural maps:
\begin{enumerate}
\item{
the evaluation map
\begin{eqnarray*}
\ovl{\cM}_{0,n}(\bP^1,1)& \longrightarrow & (\bP^1)^n; \\
(C,p_1,\ldots,p_n,\mu) & \mapsto & (\mu(p_1),\ldots\mu(p_n));
\end{eqnarray*}}
\item{
forgetting the point $p_j$
\begin{eqnarray*}
\phi_j:\ovl{\cM}_{0,n} & \longrightarrow & \ovl{\cM}_{0,n-1}, \quad (n\ge 4)\\
(C,p_1,\ldots,p_n) & \mapsto &
(C',p_1,\ldots,{\hat p_j},\ldots, p_n); \\
\phi_j:\ovl{\cM}_{0,n}(\bP^1,1) & 
\longrightarrow & \ovl{\cM}_{0,n-1}(\bP^1,1)\\
(C,p_1,\ldots,p_n,\mu) & \mapsto &
(C',p_1,\ldots,{\hat p_j},\ldots, p_n,\mu').
\end{eqnarray*}}
\item{
taking projective equivalence classes
\begin{eqnarray*}
\psi:\ovl{\cM}_{0,n}(\bP^1,1) & \longrightarrow & \ovl{\cM}_{0,n}\\
(C,p_1,\ldots,p_n,\mu) & \mapsto &
(C',p_1,\ldots, p_n).
\end{eqnarray*}}
\end{enumerate}
$C'$ is obtained from $C$ by `collapsing' the irreducible
components which are destabilized when $p_j$ (resp., the polarization) 
is removed.

Finally, we enumerate the boundary divisors of these moduli spaces.
For each partition
$$\{1,\ldots,n\}=S\cup S', \quad 2\le |S|\le |S'| \le n-2,$$
consider stable curves 
$$C=(\bP^1,p_j,j\in S) \cup (\bP^1,p_j,j\in S'),$$
which form a divisor $\delta_{S,S'}\subset \ovl{\cM}_{0,n}$.  
The union of these is denoted $\delta$.  
Note that the ${\mathfrak S}_n$-orbits of $\{\delta_{S,S'}\}$
correspond to the integers
$$|S|=2,\ldots, \lfloor n/2 \rfloor.$$
For each subset
$$S\subset \{1,\ldots,n\}, \quad 2\le S$$
consider stable maps
$$
\mu:C=(\bP^1,p_j,j\in S) \cup (\bP^1,p_j,j\in S') \longrightarrow
\bP^1
$$
collapsing the first component and mapping the second isomorphically
onto $\bP^1$.  These form a divisor $B_S\subset \ovl{\cM}_{0,n}(\bP^1,1)$.
The ${\mathfrak S}_n$-orbits of $\{B_S\}$ correspond to integers
$$
s=|S|=2,\ldots,n
$$ 
and we define 
$$
B[s]:=\sum_{|S|=s} B_S\,\,\, {\rm and}\,\,\, B:=\sum_{s=2}^n B[s]
$$

\begin{thm}
\label{thm:smooth}
The moduli spaces $\ovl{\cM}_{0,n}(\bP^1,1)$ 
and $\ovl{\cM}_{0,n}$ are smooth
projective algebraic varieties. Moreover, the boundary
is a divisor with strict normal crossings. 
\end{thm}

\begin{rem}
\label{rem:lcc}
In particular, the pair $(\ovl{\cM}_{0,n}(\bP^1,1),B)$ is log-canonical.
\end{rem}

\

\subsection{Resolution for the full moduli space}
\label{sect:res}

We obtain a good resolution $(\tilde{X},\tilde{D})$ 
of $(X,D)$ using the above formalism.
Consider the quotient map
$$
q\,:\, \ovl{\cM}_{0,n}(\bP^1,1)\ra  
\tilde{X}:=\ovl{\cM}_{0,n}(\bP^1,1)/{\mathfrak S}_n.
$$
Let $\tilde{D}[s]$ and $\tilde{D}$ be the images of $B[s]$  
and $B$ under this map.

\begin{prop}
\label{prop:ba}
The map $q$ is ramified only along the boundary $B$.
At the generic points of $B[2]$ the ramification has order $2$. 
For all $s=3,..,n$, the map $q$ is unramified 
at the generic points of $B[s]$. 
We have the formula
$$
q^*(K_{\tilde{X}} +\tilde{D}) = K_{\ovl{\cM}_{0,n}(\bP^1,1)} +B
$$
and $(\tilde{X},\tilde{D})$ is log-canonical. 
\end{prop}

\begin{proof}
The map $q$ ramifies at points corresponding to stable maps 
$$
(C, p_1,...,p_n,\mu)
$$
that admit an automorphism permuting the marked points.
The ramification order is the order of this automorphism group.
If the images of the $n$ points under $\mu$ are distinct
then there is no automorphism of $\mu$ permuting them.
This proves the first assertion. 
If marked points coincide there
is an irreducible component $\bP^1\subset C$ which is collapsed 
by $\mu$ and which contains these points. If there are two such points
this component admits an automorphism of order 
two  exchanging the points and fixing the point of intersection 
with the rest of $C$. 
This proves the second assertion.
If there are more than two
marked points then there is generally no such automorphism. 
This proves the third assertion.
 
The ramification formula and the fact that the pair
$(\tilde{X},\tilde{D})$ is log-canonical 
follow from an easy local computation combined with Remark~\ref{rem:lcc}
(see Propositions 20.2 and 20.3 of \cite{flab}).
\end{proof}

Take ${\mathfrak S}_n$-quotients of the point map to obtain
a birational map
$$
\varrho\,:\, \tilde{X} \ra \bP(\Sym^n V^*),
$$
assigning to $p_1,...,p_n\in \bP^1$ a polynomial vanishing at these
points. The boundary divisor $\tilde{D}[2]$ is the proper transform of the 
discriminant $D$ under $\varrho$. The boundary divisors 
$\tilde{D}[s]$ (for $s\ge 3$)
are the exceptional divisors for $\varrho$.

\subsection{Resolution of the generic orbit}
\label{sect:res-gen}

Let $\al:=(\al_1,...,\al_n)$ be a set of distinct complex numbers
and $f=f_{\al}$ 
the binary form of degree $n$ with roots $\al_j$.
Let $C_{\al}\in \cM_{0,n}$ be the corresponding pointed rational curve  
and $\mu_{\al}\in \cM_{0,n}(\bP^1)$ the corresponding map. 
The fiber 
$$
Y_{\al}:=\psi^{-1}(C_{\al})\subset \ovl{\cM}_{0,n}(\bP^1,1)
$$
contains $\mu_{\al}$. 
Let $\tilde{X}_f$ be the image of $Y_{\al}$ under the quotient map $q$
and $\tilde{D}_f$ its intersection with the boundary $\tilde{D}$.
This coincides with the general fiber of the map
$$
\psi'\,:\, \tilde{X}= \ovl{\cM}_{0,n}(\bP^1,1)\ra 
\ovl{\cM}_{0,n}/{\mathfrak S}_n.
$$
The map $\varrho$ induces a resolution
$$
\varrho_f\,:\, \tilde{X}_f\ra X_f,
$$ 
with $\varrho_f(\tilde{D}_f)=D_f$.

\

To describe the $Y_{\al}$ explicitly, we use the tower

\centerline{ 
\xymatrix{
\ovl{\cM}_{0,n}(\bP^1,1)\ar[d]_{\phi_n}\ar[r]^{\psi}  & \ovl{\cM}_{0,n}\ar[d]^{\phi_n}\\
\ovl{\cM}_{0,n-1}(\bP^1,1)\ar@{.}[d]\ar[r]^{\psi}& \ovl{\cM}_{0,n-1}\ar@{.}[d] \\
\ar[d]_{\phi_4}   &       \ar[d]^{\phi_4} \\
\ovl{\cM}_{0,3}(\bP^1,1)\ar[r]^{\psi}           & \ovl{\cM}_{0,3}.
}
}
\no
When $n=3$, $\ovl{\cM}_{0,3}={\rm point}$ and $Y_{\alpha}\simeq 
\ovl{\cM}_{0,3}(\bP^1,1)$, which is isomorphic to the product $(\bP^1)^3$
blown up along the small diagonal $\Delta_{\rm small}$. 
The boundary divisors correspond to the following stable maps
$$
B[2] = \begin{array}{c} \includegraphics{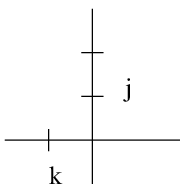} \end{array},\,\,\,\,  
B[3] = \begin{array}{c} \includegraphics{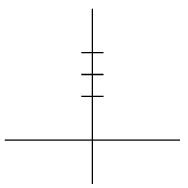}  \end{array}.
$$
In the above pictures the collapsed components are represented by 
vertical lines. 
Note that the  normal bundle
$$
\cN_{\Delta}\cong\cO(2)\oplus \cO(2),
$$
so that the exceptional divisor 
$E=B[3]\simeq \bP^1\times \bP^1$.
Let 
$$
\pi_1\,:\, E\ra \bP^1
$$ 
be the projection to 
the cross ratio of the marked points and the node 
and
$$
\pi_2\,:\, E\ra \bP^1
$$
the projection onto the image of the collapsed curve.

\

The divisor $B[2]$ is the proper transform of $\Delta$,
the large diagonal.

\

For the arbitrary degree case, we analyze the failure 
of the block squares in the tower to be fiber products.
Given a generic 
$$
C_{\al}=(\bP^1,\alpha_1,\ldots,\alpha_n)\in 
\cM_{0,n},\quad \alpha_i\ne \alpha_j,
$$
we compare the fibers 
$$
Y_{\al_1,...,\al_n}=\psi^{-1}(C_{\al}) \text{ and } 
Y_{\al_1,...,\al_{n-1}}=\psi^{-1}(\phi_n(C_{\al}))=
\psi^{-1}(\bP^1,\alpha_1,\ldots,\alpha_{n-1})
$$
using the forgetting map 
$$
\phi_n\,:\, \ovl{\cM}_{0,n}(\bP^1,1)\ra \ovl{\cM}_{0,n-1}(\bP^1,1).
$$
Given a stable map 
$$
(C,\alpha_1,\ldots,\alpha_{n},\mu) \in \psi^{-1}(C_{\al}),
$$ 
there are three cases to consider:

\begin{enumerate}
\item $C=\bP^1$;
\item $C=\bP^1\cup \bP^1$ with the collapsed component containing 
$\al_1,...,\al_n$; 
\item $C=\bP^1\cup \bP^1$ with the collapsed component containing 
$\al_1,...,\al_{n-1}$
but not $\al_n$.
\end{enumerate}

\

Case 1:   
\hskip 2cm $\begin{array}{c} 
\includegraphics{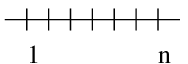} \end{array}\,\,\, \stackrel{\psi}{\longrightarrow}   \,\,\, 
\begin{array}{c} \includegraphics{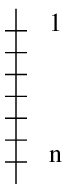} \end{array}$

\

Case 2:   
\hskip 2cm $\begin{array}{c}  
\includegraphics{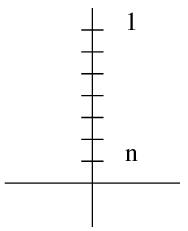} \end{array}\,\,\, \stackrel{\psi}{\longrightarrow}\,\,\, 
\begin{array}{c} \includegraphics{mo} \end{array}$

\

Case 3:   
\hskip 2cm 
$\begin{array}{c} \includegraphics{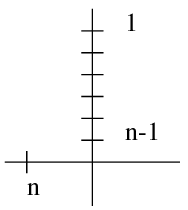} 
\end{array}\,\,\, \stackrel{\psi}{\longrightarrow}\,\,\, 
\begin{array}{c} \includegraphics{mo} \end{array}$

\

Over the open subset of $Y_{\al_1,...,\al_n}$
corresponding to the first two cases, $\phi_n$ 
induces an isomorphism between 
$Y_{\al_1,...,\al_n}$ and  
$Y_{\al_1,...,\al_{n-1}}$.  
In the third case, we forget the image
of the $n$-th marked point. 
The map $\phi_n$ blows up the locus
in $Y_{\al_1,...,\al_{n-1}}$ where $\al_1,...,\al_{n-1}$ are on the collapsed
component and $\al_n$ coincides with the node (of attachment).
This is a curve isomorphic to $\bP^1 \subset 
(B[n-1]\cap Y_{\al_1,...,\al_{n-1}})$:
The generic map takes the form:

\

\centerline{
\includegraphics{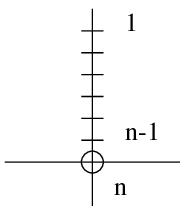}
}

\

We summarize the above discussion in the following 

\begin{prop}
\label{prop:all}
Let $\al_1,...\al_n$ be distinct complex numbers. 
The forgetting maps induce a sequence of 
birational morphisms
$$
Y_{\al_1,...,\al_n}\stackrel{ \phi_n}{\longrightarrow} Y_{\al_1,...,\al_{n-1}} \ldots 
\stackrel{ \phi_4}{\longrightarrow} 
Y_{\al_1,\al_2,\al_3}\simeq \ovl{\cM}_{0,3}(\bP^1,1).
$$
The moduli space of stable maps $\ovl{\cM}_{0,3}(\bP^1,1)$ 
is isomorphic to $(\bP^1)^3$
blown up along the small diagonal with exceptional divisor 
$E\simeq \bP^1\times \bP^1$. 
The map $\phi_j$ blows up the proper transform of 
$\pi_1^{-1}(\al_j)$.
In particular, $Y_{\al_1,...,\al_n}$ is smooth and its boundary 
has strict normal crossings, contained in $B[n-1]\cup B[n]$.
\end{prop}

\begin{rem}
\label{rem:dis}
We are blowing up along disjoint curves, so
the order of the blow-up does not matter.
\end{rem}

\begin{prop}
\label{prop:fff}
Let $f$ be a generic binary form of degree $n\ge 3$ with roots 
$\al_1,...,\al_n$. 
Then the restriction of $q$  to $X_f$ 
is ramified only along the boundary $B\cap X_f$.
At generic points of $(B[n]\cup B[n-1])\cap Y_{\al}$, 
the restriction of $q$ is unramified. 
We have the formula
$$
q^*(K_{\tilde{X}_f} +\tilde{D}_f) = K_{Y_{\al}} + [Y_{\al}\cap B]
$$
and $(\tilde{X}_f,\tilde{D}_f)$ is log-canonical. 
\end{prop}

\begin{proof}
The argument is similar to the one in Proposition~\ref{prop:ba}, and is omitted.
\end{proof}

\

\section{Verification of exponents}

\subsection{Explicit basis of $\Pic(Y_{\al})$}
\label{sect:exp}

Write
$$
\Pic((\P^1)^3)=\Z g_1+\Z g_2 +\Z g_3,\quad g_i=
{\rm pr}_i^*(c_1(\cO_{\bP^1}(+1))), 
$$
with large diagonals
$$ 
\Delta_{ij}=g_i+g_j-E,\,\, B[2]=2(g_1+g_2+g_3)-3E.
$$
By Proposition \ref{prop:all}, $Y_{\al}$ is obtained by blowing up 
the $(n-3)$ sections of 
$$
\pi_2\,:\, B[3]\ra \bP^1.
$$
Let $F_4,...,F_n$ denote the corresponding exceptional divisors and 
identify $E$ and its proper transform.
Relabel
$$
\begin{array}{ccl}
F_k & = & \Delta_{ij},\,\, \{i,j,k\}=[1,2,3]\\
    & = & g_i+g_j-E-F_4- \cdots - F_n,
\end{array}
$$
so that $\mathfrak S_n$ acts on the $F_k$, $k=1,...,n$,
in the obvious way.  Note that $E$ and the $F_k$ generate
$\Pic(Y_{\al})$.

\begin{prop}
\label{prop:ky}
The $\mathfrak S_n$-stable boundary divisors
$$
\begin{array}{rcl}
A[n-1]& = & F_1+\cdots + F_n\\
A[n] & = & E,
\end{array}
$$
generate the $\mathfrak S_n$-invariant Picard group of $Y_{\al}$,
and $A[j]=B[j]\cap Y_{\al}$.
The canonical class of $Y_{\al}$ is
$$
\begin{array}{ccl}
K  & = &  -2(g_1+g_2+g_3)+E+2(F_4+\cdots + F_n) \\
   & = & -A[n-1]-2A[n]
\end{array}
$$
\end{prop}

\subsection{Computation of the effective cone}
\label{sect:efc}

\begin{lemm}
\label{lemm:comp}
The $\mathfrak S_n$-invariant effective cone of $Y_{\al}$ is generated
by the classes $A[n]$ and $A[n-1]$. 
\end{lemm}

\begin{proof}
We apply the method of \S \ref{sect:compeff}.  
The class $A[n]$ is exceptional and thus a generator of the 
effective cone.  To show that $A[n-1]$ is the second generator,
we exhibit a nef curve not intersecting 
$A[n-1]$.  Consider the ${\mathbb G}_m$-action on $\bP^1$:
$$
\rho_t\,:\, (z,w)\mapsto (tz,w).
$$ 
We may assume that the points 
$\al_1,...,\al_n$ are not contained 
in the fixed point locus of $\rho_t$.
Any singular element in the orbit closure is:

\centerline{\includegraphics{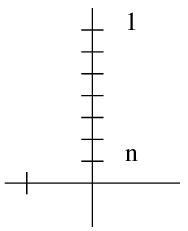}}

\no
where the point of attachment is 0 (or $\infty$) and the other 
labelled point is $\infty$ (resp. 0).  This
is disjoint from $A[n-1]$. 
\end{proof}

\subsection[Proof of the main theorem]{Proof of Theorem~\ref{thm:main}}

\begin{proof}
The Kodaira energy for 
$(X_f,D_f)$ can be computed 
on $(\tilde{X}_f,\tilde{D}_f)$, by 
Propositions~\ref{prop:ind} and \ref{prop:fff}.
By Propositions~\ref{prop:eff} and \ref{prop:all}, 
it suffices to compute the 
Kodaira energy for $(Y_{\al},A[n]+A[n-1])$.
Recall there is a composed morphism
$$
\beta:Y_{\al} \stackrel{q}{\ra} X_f \hookrightarrow 
\bP(\Sym^nV^*)\simeq \bP^n.
$$

\begin{lemm}
The pull-back of the hyperplane class takes the form
$$
L=:[\beta^*\cO_{\bP^n}(+1)]= \frac{1}{2}\left( (n-2)A[n-1]+nA[n]\right).
$$
\end{lemm}

\begin{proof}
Let $R$ be the class of a curve in $A[n]$ 
corresponding to 

\centerline{
\includegraphics{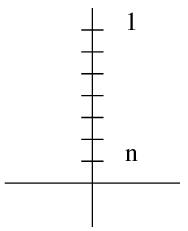}
}

\no
with varying point of attachment on the {\em collapsed} component. 
This is the proper transform of the generic fiber of the map
$\pi_1\,:\, E\ra \bP^1$. 
Then 
$$
\begin{array}{rcc}
A[n-1]|_R  & = & 
n =\# 
\left\{ \begin{array}{c} \includegraphics{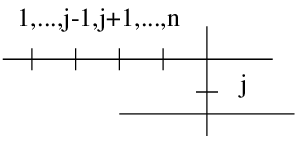} \end{array},
 		j=1,\ldots,n \right\}  \\
A[n]|_R    & = & -1-(n-3)=2-n.
\end{array}
$$
For the second intersection number, note that $A[n]=E$ and 
apply the blow-up description of Proposition \ref{prop:all}.
In $\ovl{\cM}_{0,3}(\bP^1,1)$ we have 
$$
E=\bP(\cN_{\Delta_{\rm small}})=\bP^1\times \bP^1
$$ 
and $\cN_E=\cO(-1)$.  After blowing up $(n-3)$ 
further sections of 
$$
E\ra \Delta_{\rm small}
$$
the normal bundle is reduced to $\cO(-1-(n-3))$. 

We know that $A[n]=B[n]\cap Y_{\al}$ is collapsed by the map $\beta$,
so
$$
\beta^*\cO(1)=c\left((n-2)A[n-1]+nA[n]\right)
$$
for some $c\in \N$. 
Since
\begin{multline*}
(n-2)A[n-1]+nA[n]\\
=2\left((n-2)(g_1+g_2+g_3)-(n-3)E-(n-2)(F_4+\cdots +F_n)\right).
\end{multline*}
the claim follows.
\end{proof}

We have
\begin{eqnarray*}
K_{Y_{\al}}+[Y_{\al}\cap B]&=&K_{Y_{\al}}+A[n]+A[n-1]=-A[n]\\
 K_{Y_{\al}}+[Y_{\al}\cap B]+\alpha L & =& 
\alpha \frac{n-2}{2}A[n-1] +(\alpha \frac{n}{2}-1)A[n],
\end{eqnarray*}
and by definition 
$$
a(\beta^*L):=\inf\{ \alpha \,|\, \alpha \beta^*L+K_{Y_{\al}}+[Y_{\al}\cap B]
\in \Lambda_{\rm eff}(Y_{\al})\}.
$$
Hence Lemma \ref{lemm:comp} yields
$$a(L,[Y_{\al}\cap B])=2/n.$$ 
Thus $a(L,D_f)=2/n$, as desired!
\end{proof}

\section[The invariant effective cone]{The $\fS_n$-invariant effective cone of the full moduli space}
\label{sect:fm}

In this section, we compute the ${\mathfrak S}_n$-invariant part of the 
effective cone of $\oM_{0,n}(\bP^1,1)$, 
its canonical class, and the Kodaira energy of the 
line bundle $L:=\beta^*\cO_{\bP^n}(+1)$, where
$$
\beta=\rho\circ q\,:\,  \oM_{0,n}(\bP^1,1) \stackrel{q}{\ra}
{\tilde X} \stackrel{\rho}{\ra} \bP^n.
$$
We will also compute the Kodaira energy of 
$H:=\rho^*\cO_{\bP^n}(+1)$. 

\

We first recall some basic facts about
$\oM_{0,n}(\bP^1,1)\simeq \bP^1[n]$, following \cite{fm}. 
In addition to the divisor classes $B_S$ introduced above,
we shall also consider 
$$
L_a:=\{(C,p_1,\ldots,p_n,\mu)\in \oM_{0,n}(\bP^1,1)
:\mu(p_a)=0\in \bP^1 \}, \quad a=1,\ldots,n.
$$
The cohomology $H^*(\oM_{0,n}(\bP^1,1))$ is generated by the 
classes $L_a$ and $B_S$ ,
subject to the relations
\begin{enumerate}
\item $L_a^2=0$;
\item $B_S\cdot B_{S'} =0$ for $S\cap S'\neq \emptyset$;
\item $(L_a-L_{a'})B_S=0$ for $a,a'\in S$;
\item $(\sum _{S\supset\{ a,a'\}}B_S) = L_a+L_{a'}$, 
for $1\le a < a'\le n$.
\end{enumerate}

The generators of the $\mathfrak S_n$-invariant subspace are
$$
L:=\sum_{a=1}^n L_a,\,\,\, B[s]=\sum_{|S|=s} B_S, \,\,\, 2\le s \le n.
$$
After averaging over $\mathfrak S_n$
\begin{equation}
\label{relation}
(n-1) L                     = \sum_{s=2}^n \frac{s(s-1)}{2} B[s].
\end{equation}

\begin{thm}
\label{thm:efff}
The classes $D[2], ..., D[n]$ generate the effective cone of
$Y$. 
The classes $B[2], ..., B[n]$ generate the $\mathfrak S_n$-invariant
effective cone of the moduli space  $\oM_{0,n}(\P^1,1)$.
\end{thm}

\begin{proof}
We implement the strategy of \S \ref{sect:compeff} with
$$\Gamma=\{B[2],\ldots,B[n] \}.$$
This entails finding curve classes that are 
nef relative to $\Gamma$.  
Let 
$$M=\sum_{j=2}^n d_j B[j]$$
denote an 
${\mathfrak S}_n$-invariant divisor class 
with no boundary divisors as fixed components.

Recall the description of the boundary divisor $B_S$:
$$
\begin{array}{ccll}
B_S & \simeq  & \oM_{0,n+1-s}(\P^1,1)\times \ovl{\cM}_{0,s+1},  & s=|S|>2, \\
     & \simeq & \oM_{0,n-1}(\P^1,1),                            & s=2.
\end{array}
$$
Take $s\ge 3$ and let 
$C_s\subset B_S$ be the class of the generic fiber of the map
$$
\oM_{0,n+1-s}(\P^1,1)\times \ovl{\cM}_{0,s+1}\ra 
\oM_{0,n+1-s}(\P^1,1)\times \ovl{\cM}_{0,s}
$$
forgetting the attaching point.  Since
$C_s$ passes through the generic point of $B_S$,
averaging $C_s$ over $\fS_n$ yields a curve class which
is nef relative to $\Gamma$.  In particular,
for each ${\mathfrak S}_n$-invariant divisor 
$M=\sum_{j=2}^n d_j B[j]$, moving relative to $\Gamma$, we have
$C_s\cdot M \ge 0.$

We compute intersections of $C_s$ with the various
elements of $\Gamma$.  First, the map $\beta$ blows down
the divisors $B_S$ for $|S|\neq 2$;  the data of the
collapsed component is lost completely.  It follows that
$L\cdot C_s=0$.  A simple combinatorial analysis gives
$$
C_s\cdot B_T =\left\{ \begin{array}{cl}  1      & \text{ if }
					T=S-\{ \sigma\}; \\
                                    0      & \text{ otherwise, unless } T=S. 
\end{array} \right. 
$$
which means that 
$C_s\cdot B[s-1]=s$.  Relation \ref{relation} gives
$$
0=B_S\cdot C_s \frac{s(s-1)}{2} + s \frac{(s-1)(s-2)}{2}
$$
so $B_S\cdot C_s = -(s-2)$.
To summarize, we have
$$
C_s\cdot B[j] =\left\{ \begin{array}{cl}  s      &\text{ if } j=s-1; \\
                                   -(s-2)    & \text{ if }j=s; \\
                                    0        & \text{ otherwise.} 
\end{array} \right.
$$

Using this information, we extract inequalities
on the coefficients of $M$.  The condition
$M\cdot C_s\ge 0$ yields
$$
sd_{s-1}\ge (s-2) d_s,
$$
so we get a chain of inequalities:
\begin{equation}
\label{chain}
d_n  \le  \frac{n}{n-2}d_{n-1} 
      \le  \frac{n(n-1)}{(n-2)(n-3)}d_{n-2} 
    \le  \ldots  
    \le  \frac{n(n-1)}{2}d_2.
\end{equation}
If some $d_s<0$ then $d_j<0$ for each $j\ge s$.

\

We consider another curve class in $B_S$ to
get inequalities in the reverse direction.  Fix $s\ge 2$
and let $R_s$ denote the class of the generic fiber of
$$
\oM_{0,n+1-s}(\P^1,1)\times \ovl{\cM}_{0,s+1}\ra 
\oM_{0,n-s}(\P^1,1)\times \ovl{\cM}_{0,s+1}
$$
induced by forgetting $\tau$, one of the $n+1-s$ points 
not contained in $S$.   
Again, $R_s$ passes through the generic point of $B_S$,
so averaging over $\fS_n$ yields a curve class such that
$R_s\cdot M \ge 0.$

We compute intersections as before.  The map $\beta$
sends $R_s$ to a line in $\bP^n$, 
i.e., the linear forms with $n-1$ fixed roots
and one varying root.  It follows that $L\cdot R_s=1$.  
The line $R_s$ intersects $B_T$ properly in the 
following cases
$$
R_s\cdot B_T =\left\{ \begin{array}{cl}  1      & 
	\text{ if } T=S\cup \{ \tau\}; \\
	 1     & \text{ if }
	 T=\{ \tau,\upsilon\}, \upsilon \not \in S; \\
                                    0      & \text{ otherwise, unless } T=S.
\end{array} \right. 
$$
Summing over $\fS_n$-orbits gives
$$
R_s\cdot B[j] =\left\{ \begin{array}{cl} 1  & \text{ if } j=s+1;\\    
                                        n-s-1 & \text{ if } j=2; \\
                                           0  & \text{ otherwise, unless }j=s. 
\end{array} \right.
$$
Applying Relation \ref{relation}, we find
$$
(n-1)=s(s-1)/2 R_s\cdot B[s] + (s+1)s/2 +(n-s-1),
$$
so $R_S\cdot B_S=R_s\cdot B[s] =-1.$

We extract the inequalities
$$d_{s+1}-d_s\ge (n-s-1)d_2;$$
in particular, if $d_2>0$ then each $d_j>0$.  
Adding together the inequalities
\begin{eqnarray*}
d_n-d_{n-1}&\ge & 0  \\
d_{n-1}-d_{n-2} & \ge & d_2 \\
   & \ldots & \\
d_4-d_3 & \ge & (n-4) d_2 \\
d_3 & \ge  & (n-2)d_2. 
\end{eqnarray*}
gives
$$d_n \ge \frac{n^2-5n+8}{2} d_2.$$
Combining with inequality \ref{chain}, we obtain
$$(n^2-1) d_2 \ge (n^2-5n+8) d_2,$$
hence $d_2>0$.  
\end{proof}

\begin{thm}\label{thm:full}
On the moduli space of stable maps
$(\oM_{0,n}(\bP^1,1),B)$
and its $\fS_n$-quotient $(Y,D)$, we have
$$a(B,L)=a(D,H)=2/n.$$
\end{thm}

\begin{proof}
We proceed to calculate the canonical class:
$$
K_{\oM_{0,n}(\bP^1,1)}=-2L +\sum_{s=3}^n B[s](s-2)
$$
This follows from the explicit blowup realization of 
the Fulton-MacPherson configuration space 
$\P^1[n]=\oM_{0,n}(\bP^1,1)$ \cite{fm}.
The exceptional divisors $B_S$ (for $s=|S|\ge 3$) arise from 
blowing up centers in codimension $s-1$.

\

We compute the log-Kodaira energy of $L$ on $\oM_{0,n}(\bP^1,1)$
with respect to the boundary $B$.
We have 
\begin{eqnarray*}
K+B+aL&=&(a-2)L+\sum_{s=2}^n (s-1)B[s] \\
	&=&(a-2)\left( \sum_{s=2}^n \frac{s(s-1)}{2(n-1)}B[s] \right)+
		\sum_{s=2}^n (s-1)B[s] \\
     &=& \sum_{s=2}^n (s-1) \left( \frac{s(a-2)}{2(n-1)}+1 \right) B[s]
\end{eqnarray*}  
which is effective if and only if 
$$
a\ge \frac{2(s-n+1)}{s}, \quad s=2,\ldots,n.
$$
The most restrictive inequality occurs when $s=n$,
where we obtain $a\ge 2/n$.  Consequently, 
$a(L,B)=2/n$, as expected.

The Kodaira energy for 
$(X,D)$ can be computed 
on $(\tilde{X},\tilde{D})$, by 
Propositions~\ref{prop:ind} and \ref{prop:ba}.
By Proposition~\ref{prop:eff} and Theorem~\ref{thm:smooth}, 
it is equal to the 
Kodaira energy for $(\oM_{0,n}(\bP^1,1),B)$.
\end{proof}

\section{Final remarks}

\no
A) The orbit closure $X_f$ depends on the form $[f]$. 
We get equivariant 
compactifications of $\PGL_2$ depending on moduli.
This dependence is made abundantly clear in the
blow-up description of Proposition \ref{prop:all}.

\

\no
B) The pair $(\ovl{\cM}_{0,n},\delta)$ is of log general
type:  $K_{\ovl{\cM}_{0,n}}+\delta$ is ample and log canonical
(see, for example, \S 7.1 of \cite{ha}).  
The map
$$
\psi\,:\,  \ovl{\cM}_{0,n}(\bP^1,1)\ra \ovl{\cM}_{0,n}
$$
is a log-Fano fibration onto a log-variety of general type 
(or the point, when $n=3$). 

\
 
\no
C)
The pair $(\ovl{\cM}_{0,n}, \delta) $ satisfies Vojta's conjecture.
We realize $\cM_{0,n}$ as an open subset of an algebraic
torus with explicit complement.
Fix $n-1$ points in $\bP^{n-3}$ in general position. 
Consider the set ${\mathcal H}$ of $\frac{1}{2}(n-1)(n-2)$
hyperplanes spanned by $n-3$ of the fixed points. 
Kapranov \cite{Ka} has shown that
$\cM_{0,n}\simeq \bP^{n-3} - \cup_{H\in {\cH}}H$.
The torus is obtained by excising the
$n-2$ hyperplanes spanned by subsets of the 
first $n-2$ of the points. 

\

\no
D)
We therefore expect that the asymptotic behavior of
integral points of $\ovl{\cM}_{0,n}(\bP^1,1)$ with
respect to the boundary is obtained by summing the
contributions of the integral points on the
fibers of $\psi$.  This explains
why the Kodaira energies for the moduli space of stable maps
(Theorem \ref{thm:full}) and the fibers of $\psi$
(Theorem \ref{thm:main}) should coincide 
(for the case of rational points, see \cite{BT}).  

\
\

\noindent
Brendan Hassett \\
Department of Mathematics \\
         Rice University, MS 136\\ 
	6100 S. Main Street \\
         Houston, TX 77005-1892,  U.S.A.\\
hassett@math.rice.edu

\

\noindent
Yuri Tschinkel \\
Department of Mathematics \\
         Princeton University\\ 
         Fine Hall, Washington Road\\
         Princeton, NJ 08544-1000,  U.S.A. \\
ytschink@math.princeton.edu

\end{document}